\documentclass[11pt,a4paper]{article}
\usepackage{calligra}

\usepackage[T1]{fontenc}
\usepackage[utf8]{inputenc}
\usepackage[normalem]{ulem}
\usepackage{placeins}
\usepackage{tikz}
\usepackage{amsmath}
\usepackage{amssymb}
\usepackage{graphicx}
\usepackage[amsmath,thmmarks,hyperref]{ntheorem}
\usepackage{dsfont}
\usepackage{fullpage}
\usepackage{multirow}
\usepackage[affil-it]{authblk}
\usepackage{braket}

\usepackage{enumerate}
\usepackage[shortlabels]{enumitem}

\usepackage{mathtools}
\usepackage{nicefrac}

\usepackage[backend=biber,doi=false,maxbibnames=5,style=alphabetic,isbn=false]{biblatex}

\renewbibmacro{in:}{}
\usepackage[toc,page]{appendix}
\usepackage{tocloft}
\definecolor{linkblue}{HTML}{001487}

\usepackage[colorlinks, allcolors=linkblue, breaklinks=true]{hyperref}
\theoremseparator{.}
\theorempreskip{0.7\topsep}
\newtheorem{theorem}{Theorem}

\theorembodyfont{\normalfont}
\theoremstyle{plain}
\newtheorem{definition}[theorem]{Definition}
\newtheorem*{remark}{Remark}
\makeatletter
\newtheoremstyle{MyNonumberplain}%
  {\item[\theorem@headerfont\hskip\labelsep ##1\theorem@separator]}%
  {\item[\theorem@headerfont\hskip\labelsep ##3\theorem@separator]}
\makeatother
\theoremstyle{MyNonumberplain}
\theoremheaderfont{\itshape}
\theorembodyfont{\upshape}
\theoremsymbol{\ensuremath{\Box}}

\usepackage{cleveref}

\usepackage{accents}

\DeclareMathAlphabet{\mathpzc}{OT1}{pzc}{m}{it}

\newcommand{\cP}{\mathcal{P}}

\newcommand{\bE}{\mathbb{E}}

\newcommand{\eps}{\varepsilon}

\newcommand{\R}{\mathbb{R}}

\renewcommand{\epsilon}{\varepsilon}

\renewcommand{\d}{\operatorname{d}\!}

\newcommand{\e}{\mathrm{e}}
\DeclareMathOperator{\TV}{TV}

\newcommand{\unorm}[1]{{\left\vert\kern-0.25ex\left\vert\kern-0.25ex\left\vert #1 
    \right\vert\kern-0.25ex\right\vert\kern-0.25ex\right\vert}}
\begin{document}
\title{A continuity bound for the expected number of connected components of a random graph: a model for epidemics}

\author{Koenraad Audenaert}
\affil{\small Department of Mathematics\\Royal Holloway, University of London}
\author{Eric P. Hanson}
\author{Nilanjana Datta}
\affil{\small Department of Applied Mathematics and Theoretical Physics\\University of Cambridge}
\date{December 11, 2019}

\maketitle

\begin{abstract}
We consider a stochastic network model for epidemics, based on a random graph proposed by \cite{Ros81}.
Members of a population occupy nodes of the graph, with each member being in contact with those
who occupy nodes which are connected to his or her node via edges. We prove that the expected
number of people who need to be infected initially in order for the epidemic to spread to the entire population, which is given by the expected number of connected components of the random graph,
is Lipschitz continuous in the underlying probability distribution of the random graph. We also obtain
explicit bounds on the associated Lipschitz constant. We prove this continuity bound via a technique called \emph{majorization flow} \cite{HD19a}, which provides a general way to obtain tight continuity bounds for Schur concave functions. To establish bounds on the optimal Lipschitz constant we employ properties of the Mills ratio.
\end{abstract}

Consider a graph with $n$ nodes, each representing a person, whose (undirected) edges  represent interactions which can spread infection (this is the ``two-sided epidemic process'' of \cite{Ger77}). What is the minimum number of people who need to be infected initially in order for the whole population (all $n$ nodes) to become infected eventually? The answer is simply the {\em{number of connected components}} of the graph (see Definition~\ref{compt} below). This is because if one person per connected component is infected, that person can then spread the infection to the rest of the people in the connected component. On the other hand, any connected component lacking an infected person will never become infected.

We consider the following {\em{random graph}} $G$ with $n \in \mathbb{N}$ nodes (or vertices), which was constructed in \cite{Ros81}. Fix a probability distribution $p$ on $\{1,\dotsc,n\}$, and take $n$ independent and identically distributed (i.i.d.) random variables $X_1,\dotsc,X_n \sim p$ such that $\Pr[X_i = j] = p_j$ for all $i, j \in \{1,2, \ldots, n\}$. Then construct $n$ edges by connecting $i$ to $X_i$. The result is a graph with $n$ nodes and edges, such that every node has at least one edge.

\begin{definition}\label{compt}
A connected component of a graph $G$ is a subgraph $H$ such that for every pair of nodes $x,y \in H$, there is a path (made up of contiguous edges) between $x$ and $y$, and moreover there are no edges between nodes in $H$ and $G\setminus H$ (in other words, no edges leave the connected component).
\end{definition}

\begin{figure}[ht]
\centering
\includegraphics[width=.8\textwidth]{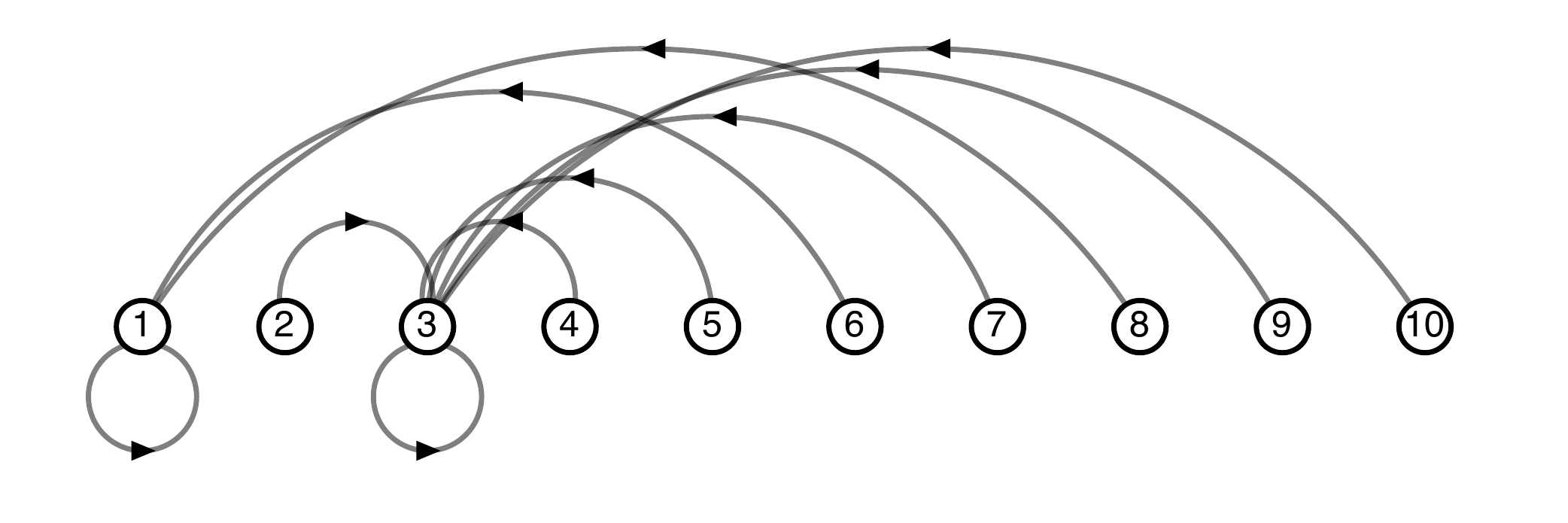}
\caption{\label{fig:onethird-twothirds} A realization of the random graph $G$ associated the distribution $p = (1/3, 0, 2/3, 0,\dotsc, 0)$. A set of edges $\{(i, X_i): i = 1,\dotsc, 10\}$ is constructed by independently sampling each $X_i \sim p$. The edges shown here are directed, but for the purpose of calculating the number of connected components, we consider the associated undirected graph.}
\end{figure}

\begin{figure}[ht]
\centering
\includegraphics[width=.8\textwidth]{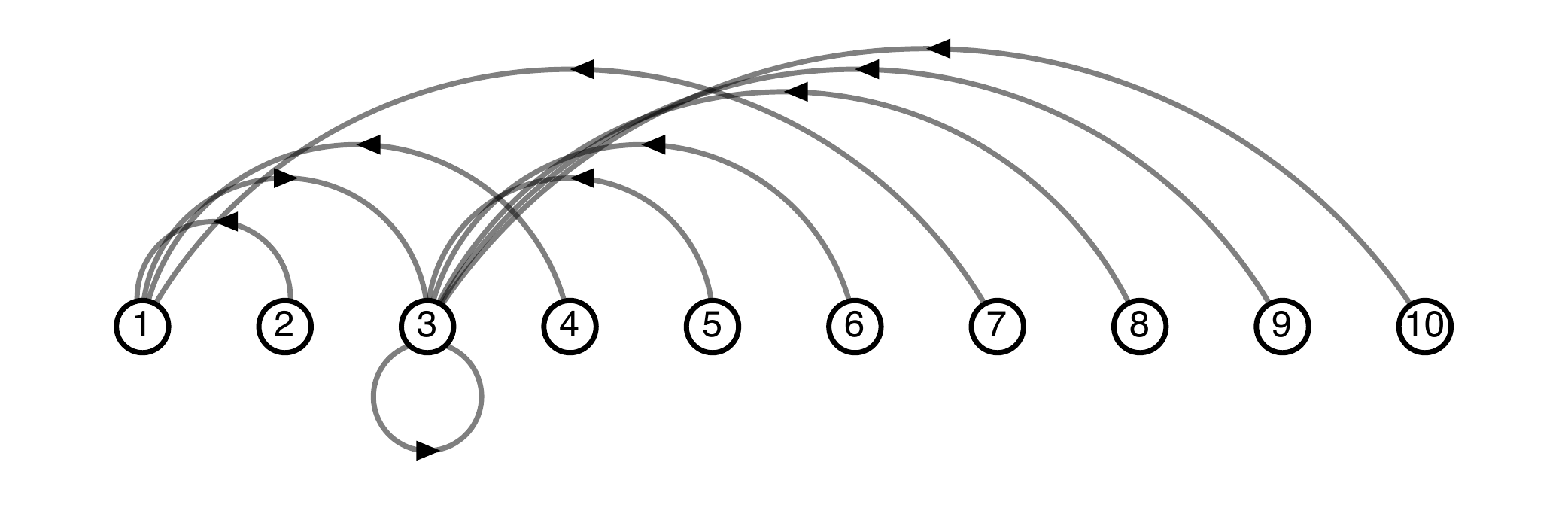}
\caption{\label{fig:onethird-twothirds-1c} Another realization of the graph $G$ from \Cref{fig:onethird-twothirds}, with underlying probability distribution $p = (1/3, 0, 2/3, 0,\dotsc, 0)$. This realization has only one connected component.}
\end{figure}

As a simple example, if $p = (1,0,\dotsc,0)$, then $1$ has a single self-edge, and every other node has a single edge touching it, which is connected to $1$, and hence the graph of one connected component. As another example, if $p = (1/3, 0, 2/3, 0,\dotsc, 0)$, then one realization of the random graph is shown in \Cref{fig:onethird-twothirds}, and another shown in \Cref{fig:onethird-twothirds-1c}.

\begin{figure}[ht]
\centering
\includegraphics[width=.8\textwidth]{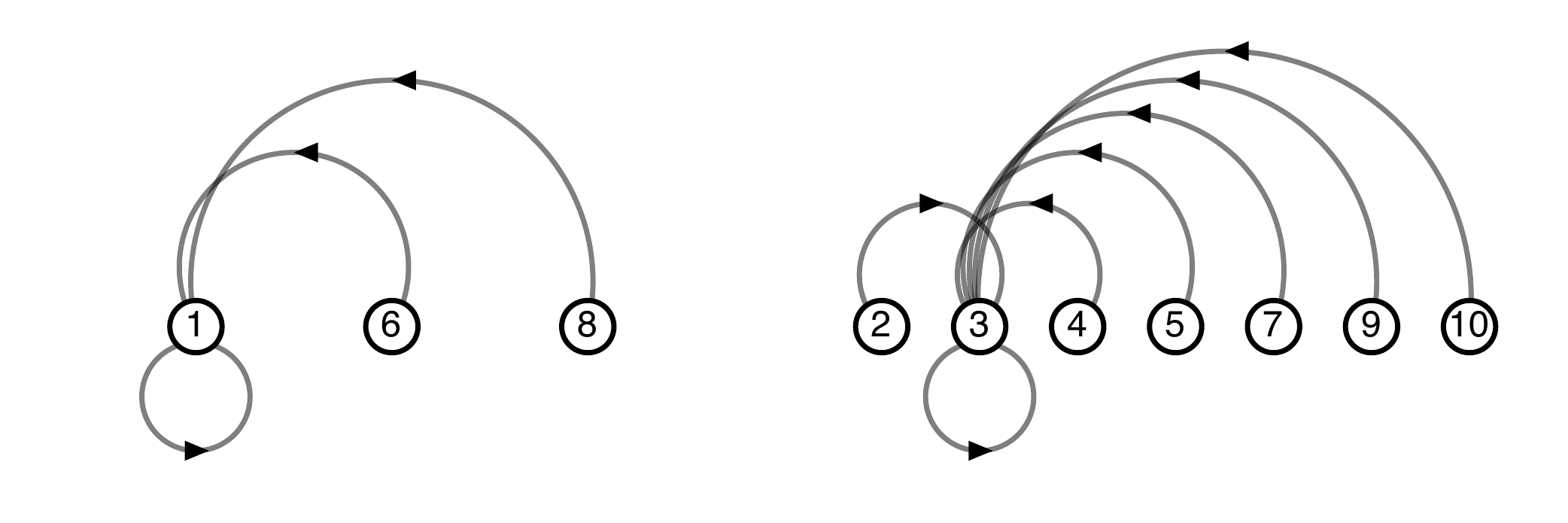}
\caption{\label{fig:onethird-twothirds-components} The two connected components of the same realization of the graph $G$ from \Cref{fig:onethird-twothirds}. These consist of elements which connect to $\{1\}$, and elements which connect to $\{3\}$, respectively. Note that we neglect directionality in computing connected components.}
\end{figure}

The number of connected components $C$ of $G$, the random graph constructed above, is a random variable. Its expectation satisfies
\begin{equation}
\bE[C] = \sum_{\substack{S\subset \{1,\dotsc, n\}\\ S \text{ nonempty}}} (|S| - 1)! \prod_{j\in S}p_j \label{eq:expected_num_connected_components}
\end{equation}
by \cite[Equation 3]{Ros81}. To make the dependence on $p$ explicit, let us write 
\begin{equation}\label{eq:def_EC}
E_C(p) \equiv \bE[C].
\end{equation}

\paragraph{Our contribution} Returning to the epidemiological framework, one can see the probability distribution $p$ as representing a model of the interactions between people, and given a probability distribution $p$, \eqref{eq:expected_num_connected_components} provides a formula for computing $E_C(p)$, the expected minimal number of initially infected people required to infect the whole population. However, if interactions are better modelled by some $p' \neq p$, then $E_C(p')$ provides a better estimate of the minimal number of people required to infect the whole population. In the following, we consider the \emph{robustness} of $E_C$; in other words, we establish a bound on the error $|E_C(p) - E_C(p')|$ as a function of the model error given by the total variation distance between $p$ and $p'$.

More specifically, we prove that the function $p \mapsto E_C(p)$ is Lipschitz continuous  on the set of probability distributions on $\{1,\dotsc,n\}$, with respect to the total variation distance, and we obtain upper and lower bounds on the optimal Lipschitz constant in terms of $n$. These bounds involve the Mills ratio,
\begin{equation}\label{eq:def_Mills_ratio}
    M(x) := \frac{1-\Phi(x)}{\phi(x)}
\end{equation}
where $\phi$ and $\Phi$ are the probability mass function and cumulative density function, respectively, of the standard normal distribution.
Our main result is given by the following theorem.

\smallskip

\noindent
\begin{theorem} \label{thm:main}
Let $E_C(p)$ and $E_C(q)$ denote the expected number of connected components of the random graph of $n \geq 3$ nodes, corresponding to probability distributions $p$ and $q$, respectively. Let $f(x) := x - x^2 M(x)$, and $\mu$ be its maximum value on the domain $x\geq 0$. Then the following bound holds:
\begin{equation}\label{eq:Lip-upper-bound}
|E_C(p) - E_C(q)| \leq  (3 + \mu \sqrt{n-2}) \TV(p, q),
    \end{equation}
where $\TV(p, q):= \frac{1}{2} ||p-q||_1\equiv \frac{1}{2} \sum_{i=1}^n |p(i) - q(i)|$ denotes the total variation distance between the two probability distributions.

Moreover, setting $x_0$ to be the unique maximizer of $f$ on $x\geq 0$, we have that any Lipschitz constant\footnote{See \eqref{eq:def_Lipschitz} below for the definition of a Lipschitz constant.} $\kappa$ satisfies
\begin{equation}\label{eq:Lip-lower-bound}
\kappa \geq \frac{\mu \sqrt{n-2} }{\sqrt{2}} - \frac{\mu x_0}{2} - \sqrt{\frac{2}{n-2}} x_0  -  \frac{x_0^2}{n-2}- \sqrt{n-2}\e^{-(n-2)}\frac{x_0\e^{x_0^2/2}}{\sqrt{2}}.
\end{equation}
Moreover, $x_0$ and $\mu$ satisfy the following explicit bounds:
\begin{align}\label{eq:bound-x0}
1.1615278892744612 &\leq x_0 \leq 1.1615278892744958,\\
0.346813047097384 &\leq \mu \leq 0.346813047097549. \label{eq:bound-mu}
\end{align}
\end{theorem}
\begin{remark}
For $n\geq 5$, the lower order terms in \eqref{eq:Lip-lower-bound} satisfy
\[
- \frac{\mu x_0}{2} - \sqrt{\frac{2}{n-2}} x_0  -  \frac{x_0^2}{n-2}- \sqrt{n-2}\e^{-(n-2)}\frac{x_0\e^{x_0^2/2}}{\sqrt{2}} \geq -2
\]
yielding the lower bound $\kappa \geq \frac{\mu \sqrt{n-2} }{\sqrt{2}} - 2$ for $n\geq 5$.
\end{remark}

\paragraph{Methods:}
The upper bound \eqref{eq:Lip-upper-bound} is proven in \Cref{sec:upper-bound}, and relies on a recent connection found between majorization, a preorder of probability distributions, and the total variation distance. This connection leads to a technique called \emph{majorization flow}, which exploits a path of probability distributions which is monotone-increasing in majorization preorder to describe the maximal change of a Schur concave function over a ball in total variation distance. This is discussed in more detail in \Cref{sec:cty-bounds}. This technique was first used by two of the present authors in \cite{HD19a} to establish novel Lipschitz continuity bounds for the $\alpha$-R\'enyi entropy, and we expect it to find use in other areas.

The lower bound \eqref{eq:Lip-lower-bound} is proven in \Cref{sec:lower-bound}. The uniqueness of the maximizer and the bounds on $x_0$ and $\mu$ are established using interval arithmetic, as described in \Cref{sec:x0mu}.

\section{Notation and definitions} \label{sec:notation}
We denote the set of probability distributions on $\{1,\dotsc,n\}$ by 
\[
\cP := \left\{ p = (p_1,\dotsc,p_n)\in\R^n : p_i \geq 0 \text{ for }i=1,\dotsc,n,\, \sum_{i=1}^n p_i = 1 \right\}.
\]
The set of probability vectors with strictly positive entries is denoted $\cP_+$. For a vector $r\in \R^n$, $r_+$ denotes its largest entry, and $r_-$ denotes its smallest entry.

A function $f: \cP \to \mathbb{R}$ is $\kappa$-\emph{Lipschitz} (with respect to the total variation distance) if for all $p,q \in \cP$,
\begin{equation}\label{eq:def_Lipschitz}
|f(p) - f(q)| \leq \kappa\, \TV(p,q)
\end{equation}
and when \eqref{eq:def_Lipschitz} holds, $\kappa$ is called a Lipschitz constant for $f$.
The smallest $\kappa>0$ such that $f$ is $\kappa$-Lipschitz is called the \emph{optimal Lipschitz constant} for $f$. The function $f$ is said to be \emph{Lipschitz continuous} if it is $\kappa$-Lipschitz for some $\kappa>0$.

Given $x\in \R^n$, write $x^\downarrow = (x^\downarrow_j)_{j=1}^n$ for the permutation of $x$ such that $x^\downarrow_1 \geq x^\downarrow_2 \geq \dotsm \geq x^\downarrow_n$. For $x,y\in \R^n$, we say $x$ \emph{majorizes} $y$, written $x \succ y$, if 
	\begin{equation} \label{def:majorize}
	 \sum_{j=1}^k x^\downarrow_j \geq \sum_{j=1}^k y^\downarrow_j \quad \forall k=1,\dotsc,n-1, \quad \text{and}\quad \sum_{j=1}^n x^\downarrow_j = \sum_{j=1}^n y^\downarrow_j.
	 \end{equation}

We say a function $\varphi: \cP \to \R$ is Schur convex if  for $p,q\in \cP$, $p\prec q\implies \varphi(p) \leq \varphi(q)$. We say $\varphi$ is \emph{Schur concave} if $-\varphi$ is Schur convex. One useful characterization of Schur convex functions is if $\varphi : \cP \to \R$ is differentiable and symmetric, then it is Schur convex if and only if
\begin{equation} \label{eq:S-convex-condition}
(p_i - p_j) \left[ \partial_{p_i}\varphi(p) - \partial_{p_j} \varphi(p) \right] \geq 0 \qquad \forall i,j
\end{equation}
for each $p \in \cP$ \cite[Section 3.A, Equation (10)]{marshall2011inequalities}.

\section{Continuity bounds for Schur concave functions}\label{sec:cty-bounds}
The quantity $E_C$ defined in \eqref{eq:def_EC} respects the majorization preorder, in that if $p, q \in \cP$ satisfy $p\prec q$, then
\[
E_C(p) \geq E_C(q)
\]
as was shown in \cite[Proposition 1]{Ros81}.
In other words, $E_C$ is Schur concave. Recently \cite{HY10,HD17,HOS18}, it has been shown that the majorization preorder interacts well with respect to total variation distance, in the sense that in any total variation ball
\[
B_\eps(r) = \left\{p \in \cP : \TV(p,r) \leq \eps \right\}
\]
there exists a minimal $r_\eps^* \in B_\eps(r)$ and maximal $r_{*,\eps} \in B_\eps(r)$ element:
\[
r_\eps^* \prec p \prec r_{*,\eps} \quad \forall p \in B_\eps(r).
\]
Moreover, in \cite{HD17}, the following so-called semigroup property was established:
\begin{equation}\label{eq:semigroup}
r_{\eps_1 + \eps_2}^* =(r_{\eps_1}^*)_{\eps_2}^* \qquad \forall r\in \cP, \, \eps_1,\eps_2 > 0.
\end{equation}
In \cite{HD19a}, this semigroup property was used to construct uniform continuity bounds for Schur concave functions with respect to the total variation distance. This construction begins by noting that if $p,q\in \cP$ satisfy $\TV(p,q) \leq \eps$, then for any Schur concave function $f$,
\begin{equation}\label{eq:bound-H-by-Delta-eps}
|f(p) - f(q)| \leq \max\left\{ f(q_\eps^*) - f(q), f(p_\eps^*) - f(p) \right\}.
\end{equation}
The semigroup property then allows the analysis of the quantity $q \mapsto f(q_\eps^*) - f(q)$ to proceed infinitesimally:
\[
f(q_\eps^*) - f(q) = \int_0^\eps \partial_s f(q_s^*) \d s = \int_0^\eps \Gamma_f(q_s^*) \d s
\]
for $\Gamma_f(r) := \left.\partial_t^+ f(r_t^*) \right|_{t=0}$, where $\partial_t^+$ denotes the derivative from above. Here, the path $(q_s^*)_{0\leq s \leq \eps}$ is the so-called path of \emph{majorization flow}. Hence, if $\Gamma_f$ is bounded above by some $k > 0$, then
\[
f(q_\eps^*) - f(q) \leq \eps k
\]
which then yields a Lipschitz continuity bound for $f$ by using \eqref{eq:bound-H-by-Delta-eps}.
Moreover, the particular structure of $r_\eps^*$ then can be used to show that for Schur concave $f$, the quantity $\Gamma_f$ is simply a difference of two partial derivatives. We refer to \cite{HD19a} for the details of this technique, which yields the following result.

\begin{theorem}[Corollary 3.2, \cite{HD19a}] \label{thm:HD19}
Let $f: \cP\to \R$ be a Schur concave function which is continuously differentiable on $\cP_+$. We write $f(r_1,\dotsc,r_n) \equiv f(r)$ for $r\in \cP$. Next, for $r \in \cP$, let $i_+ \in \{1,\dotsc,n\}$ be an index such that $r_+ = r_{i_+}$, and similarly $i_- \in \{1,\dotsc,n\}$ such that $r_- = r_{i_-}$.
Define
\begin{equation}
\begin{aligned}
\Gamma_f : \quad \cP_+ &\to \R\\
 r &\mapsto (\partial_{r_{i_+}} - \partial_{r_{i_-}})f(r_1,\dotsc, r_n).
\end{aligned}
\end{equation}
Note that this definition does not depend on the choice of $i_\pm$ since $f$ is permutation invariant. Then $f$ is Lipschitz continuous if and only if
\[
k := \sup_{r \in \cP_+} \Gamma_f(r)
\]
satisfies $k < \infty$. Moreover, in the latter case $k$ is the optimal Lipschitz constant for $f$.
\end{theorem}

\section{Proof of the upper bound \eqref{eq:Lip-upper-bound}}\label{sec:upper-bound}

In this section, we use \Cref{thm:HD19} to establish \Cref{thm:main}. Note that by \eqref{eq:expected_num_connected_components}, $E_C : \cP \to \R$ is a polynomial in the components of the probability vector $p$ and in particular is continuously differentiable.
In  \cite[Proposition 1]{Ros81}, the author proves that $p \mapsto E_C(p)$ is Schur concave using the criterion \eqref{eq:S-convex-condition}, by showing that
\[
\partial_{p_i}E_C(p) - \partial_{p_j} E_C(p) = (p_j - p_i) \sum_{S^*} (|S^*|+1)! \prod_{j\in S^*} p_j
\]
where $S^*$ ranges over nonempty sets of $\{1,\dotsc,n\}\setminus\{i,j\}$. Hence,
\[
 \Gamma_{E_C}(r) =  (r_+ - r_-) \sum_{S^*} (|S^*|+1)! \prod_{j\in S^*} r_j
\] 
where $S^*$ ranges over nonempty sets of $I :=\{1,\dotsc,n\}\setminus\{i_+,i_-\}$, where $i_\pm$ are indices such that $r_{i_\pm} = r_{\pm}$. We can use  the criterion \eqref{eq:S-convex-condition} again by repeating the proof of \cite[Proposition 1]{Ros81} to show that for
\[
S(\{r_i\}_{i\in I}) :=  \sum_{S^*} (|S^*|+1)! \prod_{j\in S^*} r_j,
\]
we have
\[
\partial_{r_i}S(\{r_i\}_{i\in I}) - \partial_{r_j} S(\{r_i\}_{i\in I}) =  (r_j - r_i) \sum_{S^*} (|S^*|+3)! \prod_{j\in S^*} r_j
\]
and hence $S$ is Schur concave on the set $\left\{ p \in \R^{n-2} : p_i \geq 0, \sum_i p_i = 1-r_+ - r_i \right\}$. For such $p$,
\[
S(p) \leq S\left( \left\{ \frac{1-r_- - r_+}{n-2} \right\}_{i \in I}\right)
\]
and thus
\begin{equation} \label{eq:EC-bound-1}
\Gamma_{E_C}(r) \leq (r_+ - r_-)\sum_{k=1}^{n-2} {n-2 \choose k} (k+1)! (1-r_- - r_+)^{k} (n-2)^{-k}.
\end{equation}
To obtain a Lipschitz bound on $E_C(r)$, it suffices to bound $\Gamma_{E_C}(r)$ independently of $r\in \cP$. We upper bound \eqref{eq:EC-bound-1} by taking $r_-=0$. For the simplicity of notation, let $s = r_+$ and $m = n-2$. Then we aim to bound
\begin{equation} \label{eq:def_Bns}
B_m(s) := s \sum_{k=1}^m{m \choose k} (k+1)! (1-s)^{k}m^{-k}
\end{equation}
for $s\in \left[\frac{1}{m+2}, 1\right]$, using that $r_+ \in \left[\frac{1}{n},1\right]$ which follows from $r \in \cP$. 
Let
$$
S_m(s) := \sum_{k=1}^m c_{k,m} (1-s)^{k-1}
$$
with
$$
c_{k,m} := {m \choose k} \frac{(k+1)!}{m^k} = (k+1)\prod_{j=1}^{k-1} \left(1-\frac{j}{m}\right),
$$
then 
\begin{equation}\label{eq:Bm-Sm-relationship}
B_m(s) = s(1-s) S_m(s).
\end{equation}

Applying the inequality $1-x\le\exp(-x)$ to every factor in $c_{k,m}$ gives the simple upper bound
\begin{eqnarray*}
c_{k,m}
&\le&  (k+1) \prod_{j=1}^{k-1} \exp(-j/m) \\
&=&(k+1)\exp\left(-\sum_{j=1}^{k-1} j/m\right)
=(k+1)\exp\left(-\frac{(k-1)k}{2m}\right) \\
&\le& (k+1)\exp\left(-\frac{(k-1)^2}{2m}\right).
\end{eqnarray*}

As $c_{k,m}\ge0$, we can also use the same inequality $1-s\le \exp(-s)$ in the formula for $S_m(s)$.
This gives as a first upper bound:
\begin{eqnarray*}
S_m(s) &\le& \sum_{k=1}^m  (k+1)\exp\left(-\frac{(k-1)^2}{2m} - (k-1)s\right) \\
&=& \sum_{l=0}^{m-1}  (l+2)\exp\left(-\frac{l^2}{2m} - ls\right) \\
&=& 2+\sum_{l=1}^{m-1}  (l+2)\exp\left(-\frac{l^2}{2m} - ls\right).
\end{eqnarray*}

We can interpret this sum as a lower Riemann sum for a certain Riemann integral.
Noting that the factor $l+2$ increases with $l$ and the factor $\exp\left(-\frac{l^2}{2m} -ls\right)$ decreases, we have
$$
(l+2)\exp\left(-\frac{l^2}{2m} -l s\right) 
\le \int_{l-1}^l (u+3) \exp\left(-\frac{u^2}{2m}-u s \right) \;du.
$$
Therefore,
\begin{eqnarray*}
S_m(s)&=&2+\sum_{l=1}^{m-1}  (l+2)\exp\left(-\frac{l^2}{2m}-ls \right)\\
&\le& 2+\sum_{l=1}^{m-1} \int_{l-1}^l (u+3) \exp\left(-\frac{u^2}{2m} -us\right) \;du \\
&=& 2+\int_{0}^{m-1} (u+3) \exp\left(-\frac{u^2}{2m}-us \right) \;du \\
&\le& 2+\int_{0}^{\infty} (u+3) \exp\left(-\frac{u^2}{2m}-us \right) \;du \\
&=& 2+\exp(ms^2/2)  \int_{0}^{\infty} (u+3) \exp\left(-\frac{(u+ms)^2}{2m} \right) \;du.
\end{eqnarray*}
In terms of the probability density function $\phi(x)$ of the standard normal distribution, $\phi(x) = \exp(-x^2/2)/\sqrt{2\pi}$,
and making the substitution $v=(u+ms)/\sqrt{m}$,
this last expression can be written as
\begin{eqnarray*}
\frac{1}{\phi(\sqrt{m}s)} \int_0^\infty (u+3) \phi\left(\frac{u+ms}{\sqrt{m}}\right)\;du
&=& \frac{\sqrt{m}}{\phi(\sqrt{m}s)} \int_{\sqrt{m}s}^\infty (\sqrt{m}v-ms+3) \phi(v)\;dv \\
&=& \frac{\sqrt{m}}{\phi(\sqrt{m}s)} \left(\sqrt{m} \int_{\sqrt{m}s}^\infty v \phi(v)\;dv + (3-ms) \int_{\sqrt{m}s}^\infty \phi(v)\;dv\right).
\end{eqnarray*}
Exploiting the fact that $x\phi(x)=-\phi'(x)$, and with $\Phi(x)$ the cumulative density function of the standard normal distribution,
this last expression is equal to
$$
\frac{\sqrt{m}}{\phi(\sqrt{m}\;s)} \left(\sqrt{m}\; \phi(\sqrt{m}\;s) + (3-mx) (1-\Phi(\sqrt{m}\;s))\right)
=m+\sqrt{m}(3-ms)\frac{1-\Phi(\sqrt{m}\;s)}{\phi(\sqrt{m}\;s)},
$$
so that
$$
S_m(s) \le 2+m+\sqrt{m}(3-ms)\frac{1-\Phi(\sqrt{m}\;s)}{\phi(\sqrt{m}\;s)}.
$$
The function in the last factor,
$$M(x):=\frac{1-\Phi(x)}{\phi(x)},$$ is known as the Mills ratio, and several bounds are known for it.
A well-known upper bound valid for $x>0$ is $M(x) < 1/x$ \cite{Gor41,YC15}, which follows from the fact that $M'(x) = x M(x) - 1$ and that $M$ is a strictly decreasing function.
Therefore, 
$$
3\sqrt{m} M(\sqrt{m}\;s)\le 3/s,
$$
and
$$
S_m(s) \le 2+m+\frac{3}{s}-m^{3/2}s M(\sqrt{m}\;s).
$$
Setting $\mu$ to be as in \Cref{thm:main}, we have
$$
-M(x)\le\frac{\mu-x}{x^2}.
$$
Therefore,
$$
S_m(s) \le 2+m+\frac{3}{s}+\sqrt{m}\frac{\mu-\sqrt{m}\;s}{s} = 2+\frac{3+\mu\sqrt{m}}{s},
$$
and
$$
B_m(s) \le (1-s)(2s+3+\mu\sqrt{m})=2(1-s)(1+s)+(1-s)(1+\mu\sqrt{m})\le 2+(1+\mu\sqrt{m}),
$$
over the interval $0\le s\le 1$

Explicit numerical calculations of $B_m(s)$ for $m$ up to $10^6$ suggest that the maximal value of $B_m(s)$ is bounded below 
by $\mu\sqrt{m}$ and, hence, lies within a constant not exceeding 3 of our bound, which is remarkable. In the following, we prove a slightly weaker bound, which recovers the square-root scaling at leading order.

\section{Proof of the lower bound \eqref{eq:Lip-lower-bound}}\label{sec:lower-bound}
Let $r = \left(r_+, \frac{1 - r_+}{n-2}, \dotsc, \frac{1 - r_+}{n-2}, 0\right) \in \cP$ for some $r_+ \in \left[\frac{1}{n-1}, 1\right]$, so that $r$ is a probability distribution with largest element $r_+$. Then the start of \Cref{sec:upper-bound} establishes that
\[
\Gamma_{E_C}(r) = B_m(s)
\]
where $B_m(s)$ is defined in \eqref{eq:def_Bns}, and $s := r_+$, and $m := n-2$. By \Cref{thm:HD19}, it remains to lower bound $B_m(s)$ for some $s \in \left[ \frac{1}{m+1}, 1 \right]$. As in \eqref{eq:Bm-Sm-relationship}, we write
\begin{equation}\label{eq:lb-decompose-B}
B_m(s) = s \sum_{k=1}^m c_{k, m} (1-s)^{k}, \qquad c_{k,m} := (k+1)\prod_{j=1}^{k-1}\left(1-\frac{j}{m}\right).
\end{equation}
Then
\[
\ln \frac{c_{k,m}}{k+1} = \sum_{j=1}^{k-1}\ln \left( 1 - \frac{j}{m} \right) = \sum_{j=0}^{k-1}\ln \left( 1 - \frac{j}{m} \right) \geq \int_0^k \ln \left( 1 - \frac{u}{m} \right) \d u
\]
using that since $j \mapsto \ln \left( 1 - \frac{j}{m} \right)$ is decreasing, the integral forms an underapproximation to the sum. By changing variables to $v = u/m$, we obtain
\[
\ln \frac{c_{k,m}}{k+1}\geq m \int_0^{k/m} \ln(1-v)\d v =  -k  - (m-k) \ln\left( 1 - \frac{k}{m} \right) \geq - \frac{k^2}{m}
\]
using that $\ln\left( 1 - \frac{k}{m} \right) \leq - \frac{k}{m}$. Hence,
\begin{equation}\label{eq:c_km-lower-bound}
c_{k,m} \geq (k+1) \exp\left( - \frac{k^2}{m}\right).
\end{equation}
From \eqref{eq:lb-decompose-B}, defining $c_{0,m} = 1$, we have
\begin{align*}
\frac{1}{s}B_m(s) &= \sum_{k=1}^{m} c_{k,m}(1-s)^k = \sum_{k=0}^{m} c_{k,m}(1-s)^k - 1\\
&\geq \sum_{k=0}^{m}(k+1) \exp\left( - \frac{k^2}{m} + k \ln(1-s) \right) - 1\\
&\geq \int_0^{m+1} u \exp \left(  - \frac{u^2}{m} + u \ln(1-s) \right) \d u - 1.
\end{align*}
using \eqref{eq:c_km-lower-bound} for the first inequality. For the second inequality, notice that the sum is of the form $\sum_{k=0}^m f(k+1) g(k)$ where $f(k)=k$ is monotone increasing, and $g(k) = \exp\left( - \frac{k^2}{m} + k \ln(1-s) \right)$ is monotone decreasing. Hence, we have $f(k+1)\geq \int_k^{k+1} f(u) \d u = \|\left.f\right|_{[k,k+1]}\|_1$ and $g(k) = \sup_{k \leq u \leq k+1} g(u) = \|\left.g\right|_{[k,k+1]}\|_\infty$, using that both functions are non-negative. H\"older's inequality gives
\[
\int_k^{k+1} f(u) g(u)\d u \leq  \|\left.f\right|_{[k,k+1]}\|_1\, \|\left.g\right|_{[k,k+1]}\|_\infty \leq f(k+1)g(k)
\]
and summing over $k$ yields the inequality. Next, since
\[
 - \frac{u^2}{m} + u \ln(1-s) = - \frac{1}{m} \left( \left( u - \frac{m \ln(1-s)}{2} \right)^2 - \left( \frac{m \ln(1-s)}{2} \right)^2\right),
\]
we obtain
\begin{align*}
\frac{1}{s}B_m(s) &\geq \frac{1}{\exp\left( - \frac{1}{2} \left(\sqrt{\frac{m}{2}}\frac{\ln(1-s)}{2}\right)^2\right)}\int_0^{m+1}u  \exp \left( - \frac{1}{2} \left( \sqrt{\frac{2}{m}}  \left(u - \frac{m}{2} \ln(1-s) \right)\right)^2\right) -1\\
&= \frac{1}{\phi(b)} \int_0^{m+1} u\phi(au -b) \d u-1
\end{align*}
for $a = \sqrt{\frac{2}{m}}$, $ b = \sqrt{\frac{m}{2}}\ln(1-s)$, and $\phi(x) := \frac{1}{\sqrt{2\pi}}\e^{-\frac{x^2}{2}}$ is the p.d.f.~of a standard normal distribution. Changing variables to $v = au - b$, we find
\begin{align*}
\frac{1}{s}B_m(s) &\geq  \frac{1}{\phi(b)a^2} \left[\int_{-b}^{a(m+1)-b} v \phi(v) \d v + b \int_{-b}^{a(m+1)-b} \phi(v) \d v\right]\\
&=  \frac{1}{\phi(b)a^2} \left[ \phi(-b) - \phi( a(m+1) - b) + b( \Phi(a(m+1)-b) - \Phi(-b))\right] - 1
\end{align*}
where $\Phi$ is the c.d.f~of the standard normal distribution. Since $a(m+1)-b \leq \sqrt{2m}$ and $\phi(x)$ is decreasing on $x > 0$, we have $\phi( a(m+1) - b) \leq \phi(\sqrt{2m})$. Using also that $\Phi( a(m+1)-b) \leq 1$, we obtain
\begin{equation*}
\frac{1}{s}B_m(s) \geq - \frac{\phi(\sqrt{2m})}{a^2 \phi(-b)} + \frac{1}{a^2} \left( 1 + b M(-b)\right) - 1 
\end{equation*}
where $M(x) = \frac{1 - \phi(x)}{\Phi(x)}$ is the Mills ratio. Substituting in $a$, we have
\begin{equation*}
\frac{1}{s}B_m(s) \geq - \frac{m\phi(\sqrt{2m})}{2\phi(-b)} + \frac{m}{2} \left( 1 + b M(-b)\right) - 1 .
\end{equation*}
Recalling the definition of $x_0$ and $\mu$ from \Cref{thm:main}, we choose $s = 1 - \e^{- \sqrt{\frac{2}{m}}x_0}$ so that $b = - x_0$, and $\mu = x_0 - x_0^2 M(x_0)$. Substituting for $\mu$, we have
\[
B_m(s) \geq s\left(\frac{m\mu}{2x_0} - 1- \frac{m\phi(\sqrt{2m})}{2\phi(x_0)} + \right) .
\]
Using the bound $\e^{-x} \leq 1 - x + \frac{x^2}{2}$ for $x\geq 0$, we have $s \geq \sqrt{\frac{2}{m}} x_0 - \frac{x_0^2}{m}$. Hence,
\[
B_m(s) \geq \frac{\mu \sqrt{m} }{\sqrt{2}} - \frac{\mu x_0}{2} - \sqrt{\frac{2}{m}} x_0  -  \frac{x_0^2}{m}- \sqrt{m}\e^{-m}\frac{x_0}{2\sqrt{\pi}\phi(x_0)}.
\]

\appendix

\section{Maximizing $x - x^2 M(x)$} \label{sec:x0mu}

The function $f(x) := x - x^2 M(x)$ on the domain $x \geq 0$ has a maximum value $\mu$ (satisfying \eqref{eq:bound-mu}) which occurs at a unique point $x_0$ (satisfying \eqref{eq:bound-x0}). To prove this, we will use the tools of \emph{interval arithmetic}. Interval arithmetic is a method for rigorous calculation using finite-precision floating point numbers on a computer, as follows. Instead of considering a real number $x\in \R$, which may not be exactly representable with a particular finite precision arithmetic, a small interval $[a,b] \subseteq \R$ containing $x$  whose endpoints are exactly representable is used instead. Then to estimate e.g. $f(x)$, instead an interval $[c,d] \subseteq \R$ is found such that $f(y) \in [c,d]$ for any $y \in [a,b]$. This yields rigorous bounds on $f(x)$ which are not subject to the ``roundoff error'' of usual floating point arithmetic. In addition, we will use the \emph{interval Newton's method}, a powerful extension of the iterative root-finding method which provides rigorous bounds on the zeros of a differentiable function and gives a sufficient condition to guarantee the function has a unique zero in a given interval \cite[Ch.~5]{Tuc11}.

First note that $f(0) = 0$ and $f(x) > 0$ using the simple upper bound $M(x) < \frac{1}{x}$ for $x>0$. Hence, any maximum of $f$ cannot occur at zero. Next, $f$ is smooth, with first derivative
\[
f'(x) = 1 + x^2 - x(2+x^2)M(x)
\]
and second derivative
\[
f''(x) =  x^3 + 4x - M(x)(2+5x^2+x^4).
\]
By using the interval Newton's method as implemented in the Julia  programming language \cite{BEKS17} package \texttt{IntervalRootFinding.jl} \cite{IntervalRootFinding.jl}, we can verify that for $x\in [0.0, 3.0]$, the equation $f'(x)=0$ has a unique solution $x_0$ which satisfies \eqref{eq:bound-x0}.
Moreover, bounding $f$ on the interval given by \eqref{eq:bound-x0} with interval arithmetic, as implemented in  \texttt{IntervalArithmetic.jl} \cite{IntervalArithmetic.jl} shows that $\mu :=f(x_0)$ satisfies  \eqref{eq:bound-mu}.
Lastly, we likewise find that
\[
-0.16730889431005824 \leq f''(x_0) \leq -0.16730889430876594,
\]
and hence $f''(x_0) < 0$ confirming that $x_0$ is a local maxima of $f$. The code used to establish these bounds can be found here: \url{https://github.com/ericphanson/AHD19_supplementary}. This code uses the MPFR library \cite{FHL+07} for a correctly-rounded implementation of the complementary error function, $1-\Phi(x)$.

Lastly, for $x > 0$, we use the lower bound $M(x) > \frac{x}{x^2+1}$ which holds for $x>0$ \cite{YC15}. This bound yields $f(x) < \frac{x}{1+x^2}$. The right-hand side is strictly monotone decreasing for $x > 1$, and evaluates to $0.3$ at $x=3$. Hence $f(x) < 0.3$ for $x > 3$. Thus, the local maximum at $x_0$ is in fact a global maximum.

\printbibliography
\end{document}